\documentstyle[11pt]{amsart}

\newtheorem{thm}{Theorem}[section]
 
\newtheorem{lem}[thm]{Lemma} 
\newtheorem{prop}[thm]{Proposition} 
\newtheorem{definition}[thm]{Definition} 
\newtheorem{rem}[thm]{Remark} 
\newtheorem{exam}[thm]{Example} 

\def\qed{\hfill $\Box$}

\newcommand\jeden {1\hskip-3.5pt1}
\newcommand\jjeden {1\hskip-2.5pt1}

\def\F{{\cal F}}

\def\Z{{\Bbb Z}}
\def\Q{{\Bbb Q}}

\def\C{{\Bbb C}}

\begin{document}

\title[Generating function of orbifold Chern classes]
{Generating Functions of Orbifold Chern Classes I:  
Symmetric Products}

\author[T.~Ohmoto]{Toru Ohmoto}
\thanks{Partially supported by Grant-in-Aid for Scientific Research (No.17340013), JSPS}
\address[T.~Ohmoto]{Department of Mathematics, 
Faculty of Science,  Hokkaido University,
Sapporo 060-0810, Japan}
\email{ohmoto@@math.sci.hokudai.ac.jp}

\maketitle

\begin{abstract}
In this paper, for a possibly singular complex variety $X$, 
generating functions of {\it total orbifold Chern homology classes} 
of the symmetric products $S^nX$ are given.  
Those are very natural ``class versions" of known 
generating function formulae 
of (generalized) orbifold Euler characteristics of $S^nX$. 
The classes work covariantly for proper morphisms.  
We state the result more generally. 
Let $G$ be a finite group and $G_n$ the wreath product $G \sim S_n$. 
For a $G$-variety $X$ and a group $A$, 
we give a``Dey-Wohlfahrt type formula" for 
{\it equivariant Chern-Schwartz-MacPherson classes} 
associated to $G_n$-representations of $A$ (Theorem 1.1 and 1.2). 
In particular, if $X$ is a point, 
this recovers a known exponential formula 
for counting  numbers $|{\rm Hom}(A,G_n)|$. 
\end{abstract}

\section{Introduction}

For a quotient variety of a  complex algebraic variety $X$ 
with an action of a finite group $G$, 
the Euler characteristic and physicists' one have been well-known: 
$$\chi(X/G)=\frac{1}{|G|} \sum_{g \in G} \chi(X^g), \quad 
\chi^{orb}(X; G)=\frac{1}{|G|} \sum_{gh=hg} \chi(X^{h,g}), $$
where $X^g$ is the set of fixed points of $g$ and $X^{h,g}:=X^h\cap X^g$, and 
the second sum runs over all pairs  $(h,g) \in G\times G$ such that $gh=hg$ 
(we may deal with non-compact varieties (but locally compact), 
so then $\chi(\cdot)$ means the Euler characteristics 
for Borel-Moore homology groups or cohomology groups with compact supports 
in analytic topology). 
Further, Bryan-Fulman \cite{BF} and Tamanoi \cite{Ta} have introduced 
a genelarization, that is, {\it orbifold Euler characteristics $\chi_m(X;G)$ 
associated to mutually commuting $m$-tuples}; $\chi=\chi_1$, $\chi^{orb}=\chi_2$. 
The aim of this paper is to generalize those kinds of characteristic numbers to 
certain ``equivariant Chern classes" \cite{O}  in connection with 
{\it classical enumerative problems in group theory} (\cite{W}, \cite{Y}, \cite{Y2}). 
In particular, we focus on the case of symmetric products. 
We also mention a bit about an algebro-geometric aspect relating to crepant resolutions 
(Remark \ref{stringy}).

To begin with, let us recall Chern homology classes of singular varieties $X$ 
(MacPherson \cite{Mac}, Schwartz \cite{Sch}). 
Let $\F(X)$ be the abelian group of constructible functions over $X$ and 
$H_*(X)$  the Borel-Moore  homology group: $\F$ and $H_*$ are covariant functors. 
It is proved in \cite{Mac} that 
there is a unique {\it natural transformation} $C_*: \F(X) \to H_*(X)$ 
such that if $X$ is nonsingular, then $C_*(\jeden_X)=c(TX) \frown [X]$. 
Here $\jeden_X$ is the constant function $1$ over $X$ and $c(TX)$ is 
the total Chern (cohomology) class of the tangent bundle. 
The total homology class $C_*(X):=C_*(\jeden_X)$ is called 
 {\it the Chern-Schwartz-MacPherson class} of $X$. 
In particular,  its $0$-dimensinal component  (the degree) coincides 
with $\chi(X)$ for compact varieties $X$. 
We remark that the Chern class theory is available also in 
a purely algebraic context over a field of characteristic $0$ (Kennedy \cite{Ken}), 
where the homology is replaced by the Chow group $A_*(X)$.  

Based on the above theory $C_*$, 
the author \cite{O} introduced 
{\it the equivariant Chern-MacPherson transformation} $C^G_*$ 
for possibly singular varieties with actions of an algebraic group $G$, 
that will be reviewed in \S 2.

A particular interest arises in 
the $n$-th symmetric product $S^nX$ of a possibly singular complex variety $X$, 
that is, 
the quotient of $X^n=X\times \cdots \times X$ ($n$ times)  via 
the $n$-th symmetric group $S_n$ with the action permuting the factors. 
There is a well-known formula due to Macdonald \cite{M}: 
\begin{equation}\label{Macdonald}
\sum_{n=0}^\infty \chi(S^nX) z^n = (1-z)^{-\chi(X)}. 
\end{equation}
We realize the ``Chern class version" of the formula (\ref{Macdonald})   
in $\sum_{n=0}^\infty z^n H_*(S^nX;\Q)$ 
(formal power series whose coefficients are total homology classes), 
which becomes a commutative and associative graded  ring 
with the cross product multiplication  $\odot$ (see \S3): 
\begin{equation}\label{MacdonaldChern}
\sum_{n=0}^\infty C_*(S^nX) z^n 
= (1-zD)^{-C_*(X)}. 
\end{equation}
Here $D$ is the notation indicating {\it diagonal operators}: 
Its ``$n$-th power" $D^n$  means the  homomorphism 
of homology groups  induced by the diagonal embedding 
${\mit\Delta}^n: X \to {\mit\Delta} X^n \subset X^n$. 
As our convention (see Subsection \ref{Diagonal}), 
we let $(1-zD)^{-c}$ denote 
$$
exp(-Log(1-zD)(c)) 
= exp(zD(c)) \odot exp\left(\frac{z^2D^2(c)}{2}\right) \odot \cdots .
$$
The ``$0$-dimensional part" of (\ref{MacdonaldChern}), 
that is, the power series whose coefficients are the $0$-dimensional parts of the total classes, 
gives (\ref{Macdonald}). 
The formula next to (\ref{MacdonaldChern}) is 
\begin{equation}\label{HHChern}
\sum_{n=0}^\infty C_*^{orb}(S^nX) z^n 
= \prod_{k=1}^\infty (1-z^kD^k)^{-C_*(X)}  
\end{equation}
whose  $0$-dimensional part coincides with 
the known generating function of $\chi^{orb}(X^n; S_n)$ (cf. Hirzebruch-H\"ofer \cite{HH}). 
Here $C^{orb}_*$ means 
our simplest orbifold Chern class (Example \ref{mthcanonical}), whose degree is just $\chi^{orb}$. 
In fact, more generally, given a group $A$  we introduce  
{\it the canonical orbifold Chern classes 
associated to group representations of $A$}  (Definition \ref{canonicalconst}).  
It is defined by the (equivariant) $C_*$-image of 
a `canonical' constructible function assigning 
to each point  (as of the quotient stack) 
the number of representations of $A$ into its stabilizer group. 
 For the symmetric product case, the class is denoted by  
 $C^{S_n}_*(\jeden_{X^n/S_n}^{(A)})$. 
We then show  the following ``Dey-Wohlfahrt formula", 
which gives the above (\ref{MacdonaldChern}) and  (\ref{HHChern}) 
as typical examples when $A=\Z$ and $\Z^2$ respectively.  
For any positive integer $r$, let $\Omega_A(r)$ (resp. $\Omega_A$) denote 
the set of all subgroups $B$ of index $|A:B|=r$ 
(resp. subgroups of finite index) and 
$j_r(A):=|\Omega_A(r)|$, the number of the subgroups of index $r$. 
\begin{thm} \label{DW}
Assume that $j_r(A) < \infty$ for any $r$. 
Then it holds that 
$$\sum_{n=0}^\infty C^{S_n}_*(\jeden_{X^n/S_n}^{(A)}) z^n 
= exp\left(\sum_{B \in \Omega_A} 
\frac{1}{|A:B|} (zD)^{|A:B|} C_*(\jeden_X)\right).$$
\end{thm} 
When $X=pt$, this theorem gives the enumerative formula for numbers 
$\frac{1}{|S_n|}|{\rm Hom}(A,S_n)|$ 
(\cite{W}; \cite{Sta}, Prob. 5.13, \cite{Y}).   
In case that $A=\Z^m$, this extends 
the generating functions for $\chi_m(X^n;S_n)$ 
given in \cite{BF} and  \cite{Ta}
(Examples \ref{freeabelian1}, \ref{freeabelian2}).

There is a certain  $G$-version of the above theorem. 
For simplicity, we assume that $G$ is a finite group, and 
let $X$ be a $G$-variety. 
We denote by $G_n(=G \sim S_n)$ {\it the wreath product}. 
$G_n$ acts on $X^n$ in an obvious way. 
\begin{thm} \label{DWG} 
Let $G$ be a finite group, and $A$ a group 
so that $j_r(A)<\infty$ and 
$|{\rm Hom}(B,G)|<\infty$ for any subgroup $B$ of finite index. 
Then it holds that 
$$\sum_{n=0}^\infty C^{G_n}_*(\jeden_{X^n/G_n}^{(A)}) z^n 
= exp\left(\sum_{B \in \Omega_A} 
\frac{1}{|A:B|} (zD)^{|A:B|} C^G_*(\jeden_{X/G}^{(B)})\right).$$
\end{thm} 
Obviously, when $G=\{e\}$, this coincides with Theorem \ref{DW}. 
Taking $X=pt$ again, this formula coincides 
with an exponential formula (M\"uller \cite{Mu}) 
for numbers of $G_n$-representations (Remark \ref{Muller}).

We end this introduction by a few remarks (see Remark \ref{stringy} below). 
As  known,  Batyrev \cite{Baty}  showed that 
$\chi^{orb}(X;G)$ coincides with 
his {\it stringy Euler number} of the quotient variety 
and thus it equals 
the Euler characteristics of any crepant resolution, 
e.g., $\chi^{orb}(X^n, S_n)=\chi({\rm Hilb}^n X)$ for a smooth surface $X$,  \cite{HH}. 
An analogue  for $C^{orb}_*$ holds 
by passing through {\it the stringy Chern class} in 
de Fernex-Lupercio-Nevins-Uribe \cite{FLNU} and  Aluffi \cite{Aluffi}, 
and therefore, for instance,  
the above formula (\ref{HHChern})  relates to 
a recent work of Boissi\`ere \cite{Bois} 
on Chern classes of ${\rm Hilb}^n X$  for a smooth surface. 
 
Futhermore,  also related with stringy class invariants, 
 a unified theory of (additive) characteristic classes for singular varieties, 
 {\it the Hirzebruch class} (and {\it the motivic Chern class}), has appeared in  
Brasselet-Sh\"urmann-Yokura \cite{BSY} (\cite{SY} for a survey). 
It is a natural transformation $T_{y*}$  from 
 {\it the relative Grothendieck ring} $K_0({\rm Var}_\C/X)$ 
 (or $\hat{\cal M}({\rm Var}_\C/X)$)  to $H_*(X)\otimes \Q[y]$, 
 which unifies Chern-MacPherson's transform $C_*$ ($y=0$), 
 Baum-Fulton-MacPherson's Todd class transform ($y=-1$)  and 
 Cappell-Shaneson's $L$-class transform ($y=1$). 
Probably, based on our construction this transformation 
may produce  variants of above type generating function formulas for 
 other characteristic classes.

\section{Review on equivariant Chern classes} 

\subsection{Natural transformation} 
Let $G$ be a complex linear algebraic group, and 
$X$ a complex algebraic variety with an algebraic action of $G$. 
We always assume that $X$ is $G$-embeddable, that is,  it admits 
a closed equivariant embedding into a $G$-nonsingular variety  
(working in the context over arbitrary base field $k$ of characteristic $0$, 
we assume $X$ is a quasi-projective variety with a linearlized action). 
As a remark, the ``quotient" $X/G$ is no longer a variety in general 
but is a {\it quotient stack}. 

A constructible function over $X$ is 
a function $\alpha: X \to \Z$ which is written  (uniquely) 
as a finite sum 
$\alpha = \sum a_i\jeden_{W_i}$ for 
some subvarieties  $W_i$ in $X$ and $a_i \in \Z$, 
where $\jeden_{W}$ means 
the characteristic function taking value $1$ over $W$ and $0$ otherwise 
(After \S 3, we will consider only the case of rational coefficients, $a_i \in \Q$).  
{\it The integral} of  $\alpha$ over $X$ is given by  
$\int_X \alpha := \sum a_i \chi(W_i)$. 
Let $\F(X)$ be the group of constructible functions over $X$ and 
$\F^G_{inv}(X)$ the subgroup consisting of 
all $G$-invariant constructible functions 
($\alpha(g.x)=\alpha(x)$ for any $g \in G$, $x\in X$). 
By definition any invariant function is uniquely expressed by a linear combination 
$\jeden_{W_i}$ for some invariant reduced schemes $W_i$.  
For a proper $G$-equivariant morphism $f: X \to Y$, 
we define {\it the pushforward} $f_*:\F^G_{inv}(X) \to \F^G_{inv}(Y)$ by 
$f_*(\alpha)(y) := \sum a_i \, \chi(W_i \cap f^{-1}(y))$ for $y \in Y$. 
Given another proper $G$-morphism $g: Y \to Z$, $g_*\circ f_*=(g\circ f)_*$ holds. 

{\it The $G$-equivariant (Chow) homology group}  $H^G_*(X)$ 
($A^G_*(X)$) is defined by Totaro \cite{Totaro}, Edidin-Graham \cite{EG} 
using an algebraic version of the Borel construction. 
It satisfies naturally expected properties, for instance,  it admits 
the equivariant fundamental class $[X]_G \in H^G_{2n}(X)$ ($n=\dim X$) 
so that 
$\frown [X]_G: H_G^{i}(X) \to H^G_{2n-i}(X)$ is isomorphic if $X$ is nonsingular. 
Throughout, we deal only with  $G$-invariant cycles in $X$ (i.e., of non-negative dimension), 
not  general ``$G$-equivariant" cycles. 

Both of $\F^G_{inv}$ and $H^G_*$ become covariant functors  
for the category of complex $G$-varieties and proper $G$-morphisms. 
Here we state the main theorem in \cite{O} but in a bit weaker form: 
\begin{thm}\label{C} 
{\rm (Equivariant MacPherson's transformation \cite{O})} 
There is a natural transformation between these covariant funtors 
$$C_*^G: \F^G_{inv}(X) \to H^G_*(X)$$
so that  if $X$ is non-singular, then $C_*^G(\jeden_X)=c^G(TX) \frown [X]_G$ 
where $c^G(TX)$ is $G$-equivariant total Chern class of the tangent bundle
of $X$. 
\end{thm}

The following elementary properties are easily checked: 

\begin{enumerate}
\item[(i)] 
{\it Trivial action}: 
If $G$ acts trivially on $X$, $C^G_*$ coincides with $C_*$ \cite{Mac}. 
\item[(ii)] 
 {\it Degree}: 
When $X$ is projective and irreducible, taking  the pointed map $pt:X \to \{pt\}$, 
{\it the degree of $C_*(\alpha)$} is defined to be $pt_*C^G_*(\alpha) \in H^G_0(pt)=\Z$,  
that equals $\int_X \alpha$. 
In particular,   the degree of $C^G_*(\jeden_X)$ is $\chi(X)$ and its top component is $[X]_G$. 
\item[(iii)] 
{\it Quotient}: 
There is a homomorphism $i^*: H^G_*(X)  \to H_*(X)$  
induced by the inclusion $i$ of $X$ 
onto a fibre of the universal bundle $X\times_G EG \to BG$ (see \cite{EG}, \cite{O}).  
If $X/G$ is a variety, 
the pushforward induced by the projection $\pi:X \to X/G$ makes sense 
(we denote also by $\pi_*$ the composition $\F^G_{inv}(X) \subset \F(X) \to \F(X/G)$), 
and then the diagram commutes 
$$
\begin{array}{ccc}
 \F^G_{inv}(X) & \stackrel{C^G_*}{\longrightarrow} & H^G_*(X)  \\
 \pi_*\, \downarrow \;\;\;   & & \; \;  \downarrow \,  \pi_*\circ i^*   \\
 \F(X/G) & \stackrel{C_*}{\longrightarrow} & H_*(X/G)  
\end{array}
$$
If $G$ is a finite group,  
the right-sided vertical map $\pi_*\circ i^* $ 
 is an isomorphism within rational coefficients $\Q$ 
(see  \cite{EG}, Thm. 3).  
\item[(iv)]  
{\it Change of groups}: 
Let $H$ be a subgroup of $G$ with $\dim G/H=k$. 
Then, natural isomorphisms $\phi_F$ and $\phi_H$ are defined in an obvious way 
so that the following diagram commutes: 
$$
\begin{array}{ccc}
\F_{inv}^{H}(X) & \stackrel{C^{H}_*}{\longrightarrow} & H^{H}_{2n-*}(X) \\
 \phi_F \; \simeq\downarrow \;\; & & \;\; \downarrow \simeq \; \phi_H \\
\F_{inv}^{G}(X\times_{H} G)  & \stackrel{C^{G}_*}{\longrightarrow} 
& H^{G}_{2(n+k)-*}(X\times_{H} G) 
\end{array}
$$
\item[(v)] 
{\it Cross product}: 
The equivariant homology admits {\it the cross (exterior) product} (\cite{EG}). 
The exterior product of equivariant constructible functions $\alpha$ and $\beta$ of $X$ and $Y$, respectively, 
is given by 
$\alpha \times \beta(x,y) := \alpha(x) \cdot \beta(y)$. 
It holds as same as the ordinary case (\cite{Kw}) that 
$$
C^G_*(\alpha \times \beta) = C^G_*(\alpha) \times C^G_*(\beta).
$$ 
\end{enumerate}

\subsection{Canonical Chern classes} \label{Canonical} 
For simplicity we assume that $G$ is a finite group. 
Let $A$ be a group so that ${\rm Hom}(A,G)$ is a finite set. 
For any  $G$-representation $\rho \in {\rm Hom}(A,G)$, 
we set $X^{\rho(A)}:=\bigcap_{g\in \rho(A)} X^{g}$ in a set-theoretic sense, 
more precisely, $X^{\rho(A)}$ is 
the reduced scheme (or the underlying reduced analytic space) 
of the fixed point set of the action of $\rho(A)$ on $X$. 
\begin{definition} \label{canonicalconst} 
{\rm We define  
{\it the canonical constructible functions of a $G$-variety $X$ 
associated to a group $A$} by 
$$\jeden_{X/G}^{(A)}=\frac{1}{|G|} \sum_{\rho}  \jeden_{X^{\rho(A)}} 
\; \in \; \F^G_{inv}(X)\otimes \Q$$
the sum being taken over all $\rho \in {\rm Hom}(A,G)$. 
We call the class 
$C^G_*(\jeden_{X/G}^{(A)}) \in H_*^G(X;\Q)$ 
{\it the canonical quotient (orbifold) Chern classes of $X$ associated to $A$} 
(we use the word ``orbifold" when $X/G$ is a variety). 
}
\end{definition} 
For $x \in X$, the value $\jeden_{X/G}^{(A)}(x)$ equals 
$\frac{1}{|G|} |{\rm Hom}(A, Stab_G(x))|$  
(hence its $G$-invariance is clear).  
So the canonical constructible function measures by using a fixed group $A$ about 
how `large' each automorphism group is.

\begin{exam}\label{mthcanonical}
{\rm 
In a typical case that $A=\Z^m$, 
we simply denote the associated function by $\jeden_{X/G}^{(m)}$, 
called {\it the $m$-th canonical function}. 

Since a representation $\Z^m \to G$ uniquely corresponds to 
a mutually commuting $m$-tuple $(g_1, \cdots, g_m)$ of $G$, 
the integral over $X$ (=degree of its canonical Chern class) coincides with the definition of 
the orbifold Euler characteristics given in \cite{BF}, \cite{Ta}: 
$$
\int_X \jeden_{X/G}^{(m)} 
=\frac{1}{|G|} \sum \chi(X^{g_1, \cdots, g_m}) =: \chi_m(X;G).
$$
Let ${\cal X}:=X/G$ be a variety. As for pushforward  via $\pi: X \to {\cal X}$, 
it is straightforward that 
$\pi_*(\jeden_{X/G}^{(1)})=\jeden_{\cal X}$: 
for any $[x] \in {\cal X}$ 
$$
\pi_*(\jeden_{X/G}^{(1)}) ([x])=\int_{G.x} \jeden_{X/G}^{(1)} 
= |G.x|\frac{|{\rm Hom}(\Z, Stab_G(x))|}{|G|} =1.
$$
Hence $C^G_*(\jeden_{X/G}^{(1)})$ is identified  
with the ordinary Chern-SM class $C_*({\cal X})$ of the quotient variety 
through $\pi_*i^*$ in property (iii) within rational coefficients. 
Furthermore,   we put  
$$
C_*^{orb}({\cal X}) :=\pi_*i^*C^G_*(\jeden_{X/G}^{(2)}) \; \in \; H_*({\cal X};\Q). 
$$
We decompose  $\pi_*\jeden_{X/G}^{(2)}$  
in the exactly same way as \cite{HH} (\cite{FLNU}) and obtain an alternative expression 
$$
C_*^{orb}({\cal X}) = \sum (\iota_g)_*C_*(X^g/C(g)),
$$
where 
the sum runs over the set of all conjugacy classes in $G$, 
$g$ is a representative in each conjugacy class, 
$C(g)$ is the centralizer of $g$, 
and $\iota_g: X^g/C(g) \to X/G$ is the canonical inclusion. 
}
\end{exam}

\begin{rem} \label{stringy}
{\rm 
We give a short remark about a connection with resolutions. 
A further account will be discussed somewhere else. 

For normal varieties with ``tame" singularities 
the {\it stringy Chern class} $c_{str}$ has been introduced by 
de Fernex et al \cite{FLNU} and  Aluffi \cite{Aluffi} (also see \cite{BSY}).  
Roughly, 
it is defined by the $C_*$-image of a constructible function 
coming from a relative motivic intergration associated to resolutions of singularities. 
Assume that the quotient variety ${\cal X}=X/G$ ($X$ being smooth) 
adimits a crepant resolution $f: Y \to {\cal X}$. 
It then follows from \cite{FLNU} (Thm. 0.2 and Thm.4.4) that  
$$
f_*(C_*(\jeden_Y))=c_{str}({\cal X})=C_*^{orb}({\cal X})  
\quad (\mbox{i.e.,} \;  f_*(\jeden_Y)=\pi_*\jeden_{X/G}^{(2)}).
$$
A naive question is to ask if 
our orbifold Chern class associated to a group $A$ has a similar property, e.g., 
if $\pi_*\jeden_{X/G}^{(A)}$ coincides with the $f_*$-image of some 
some distinguished constructible funtion on $Y$.

Let $X$ be a smooth surface. 
It is well known that the Hilbert scheme ${\rm Hilb}^n(X)$ is smooth 
and the Hilbert-Chow morphism $\pi^n: {\rm Hilb}^n(X) \to S^nX$ becomes a crepant resolution. 
Recently, 
a generating function of total Chern classes  of the tangent bundle $T\, {\rm Hilb}^n(X)$ 
was given by Boissi\`ere \cite{Bois} (Prop. 3.12) using vertex algebras tools in \cite{Lehn}. 
Since $(\pi^n)_*(C_*({\rm Hilb}^n(X)))=C_*^{orb}(S^nX)$ just as mentioned, 
the formula (\ref{HHChern}) in Introduction must be the image  
of the generating function formula via $(\pi^n)_*$ ($n\ge 0$) (and the Poincar\'e dual). 
This may suggest  a direct connection betwen our formulae  and  
exponential formulae in vertex algebras. 
From our equivariant viewpoint, 
Fulton-MacPherson compactifications of configuration spaces (which has $S_n$-actions)  
should be interesting.

}
\end{rem}

\section{Symmetric products} 

\subsection{Formal power series} \label{Csym}
We work on symmetric products, that is the quotient via 
the action of $S_n$ on the cartesian product $X^n$ 
of a complex variety $X$, 
$\sigma(x_1, \cdots, x_n)
:=(x_{\sigma^{-1}(1)}, \cdots, x_{\sigma^{-1}(n)})$. 
From now on, we deal with 
$\F$ and $H_*$ 
with rational coefficients 
and omit the notation $\otimes \, \Q$. 

The product 
$\odot: \F^{S_m}_{inv}(X^m) \times \F^{S_n}_{inv}(X^n)
\to \F^{S_{m+n}}_{inv}(X^{m+n})$ is defined to be 
$$
\alpha\odot \beta
:= \frac{1}{|S_{m+n}|} \sum_{\sigma\in S_{m+n}} \sigma_*(\alpha\times \beta), 
$$
where $\sigma_*$ is the pushforward induced 
by $\sigma: X^{m+n} \to X^{m+n}$. 
Also for homologies,  $\odot$ is defined in the same manner. 
This yields commutative and associative graded $\Q$-algebras of formal power series 
$$
\F_{X, sym}[[z]]:=\sum_{n=0}^\infty z^n \F^{S_n}_{inv}(X^n), \quad 
H_{X, sym}[[z]] :=\sum_{n=0}^\infty z^n H^{S_n}_{2*}(X^n).
$$
We denote $\alpha \odot \cdots \odot \alpha$ ($c$ times) 
by $\alpha^c$ or $\alpha^{\odot c}$. 

For a proper morphism $f:X\to Y$, 
the $n$-th cartesian product $ f^n:X^n \to Y^n$ 
 is a $S_n$-equivariant map, and it hence induces 
$$
f^{sym}_*:\F_{X, sym}[[z]] \to \F_{Y,sym}[[z]], \quad
f^{sym}_*(\sum_{n=0}^\infty \alpha_n z^n)
:=\sum_{n=0}^\infty  f^n_*\alpha_n z^n.
$$
as well the homology case. 
It is easy to see that $f^{sym}_*$ is a homomorphism of algebras 
(it preserves the multiplication).  

By definition, when $X=pt$, $\F_{pt, sym}[[z]]$ is 
canonically isomorphic to $\Q[[z]]$, 
the ring of formal power series with rational coefficients. 
{\it The integral} of a power series of constructible functions 
is defined by 
$$
\int: \F_{X, sym}[[z]] \to \Q[[z]], \quad 
\sum_{n=0}^\infty  \alpha_n z^n 
\mapsto \sum_{n=0}^\infty  (\int_{X^n} \alpha_n)\; z^n.
$$
We also define 
$$
C^{sym}_*: \F_{X, sym}[[z]]\to H_{X, sym}[[z]], \quad 
\sum_{n=0}^\infty \alpha_n z^n 
\mapsto \sum_{n=0}^\infty  C^{S_n}_*(\alpha_n) z^n.
$$
\begin{thm} \label{SymC}
$C^{sym}_*$ is a natural transformation 
between the covariant functors assigning to $X$ 
the $\Q$-algebras $\F_{X, sym}[[z]]$ and $H_{X, sym}[[z]]$. 
\end{thm}
\noindent{\it Proof}.  
We show that 
\begin{equation}\label{pc}
C^{S_m}_*(\alpha) \odot C^{S_n}_*(\beta) 
= C^{S_{m+n}}_*(\alpha\odot \beta).
\end{equation}
Put $G=S_{m+n}$ and $H=S_m\times S_n$. 
$H$ acts on $X^m$, $X^n$, $X^{m+n}$ and (the left action) on $G$, then by $(iv)$, 
$\F^{H}_{inv}(X^m)=\F^{S_m}_{inv}(X^m)$ and 
$\phi_F: \F^{H}_{inv}(X^{m+n}) \simeq 
\F^{G}_{inv}(X^{m+n}\times_{H}G)$. 
Let $p: X^{m+n}\times_{H}G \to X^{m+n}$ 
 be the natural projection given by $p([x,a]_{H}):= a^{-1}.x$, 
which is well-defined and $G$-equivariant 
(the action of the mixed space is given by $g.[x,a]_{H} :=[x,ag^{-1}]_{H}$). 
We denote by $\tau_F$  the following composed homomorphism 
\begin{eqnarray*}
&&
\F^{H}_{inv}(X^m)\otimes \F^{H}_{inv}(X^n) 
\stackrel{\times}{\longrightarrow} 
\F^{H}_{inv}(X^{m+n}) \\
&&\qquad \quad \qquad \stackrel{\phi_F}{\longrightarrow} 
\F^{G}_{inv}(X^{m+n}\times_{H}G) 
\stackrel{p_*}{\longrightarrow} 
\F^{G}_{inv}(X^{m+n}).
\end{eqnarray*}
For homology we have $\tau_H$ in the same way. 
Those $\tau_F$ and $\tau_H$ actually coincide with 
the $\odot$-products up to a scalar multiple: 
%
%
in fact, for $x \in X^{m+n}$, 
\begin{eqnarray*}
\tau_F(\alpha, \beta)(x)
&=& \int_{p^{-1}(x)} \phi_F(\alpha \times \beta) 
=\sum_{[g] \in H\backslash G} \phi_F(\alpha \times \beta)([g.x,g]_{H}) \\
&=&  \frac{1}{|H|} \sum_{g \in G}   (\alpha \times \beta)(g.x)
= \frac{(m+n)!}{m!n!} (\alpha \odot \beta) (x), 
\end{eqnarray*}
the homology case as well. 
By properties  (v), (iv) and the naturality, we see 
\begin{eqnarray*}
&&\tau_H(C_*^{H}(\alpha), C_*^{H}(\beta))
= p_*\circ \phi_H(C_*^{H}(\alpha)\times C_*^{H}(\beta)) 
 = p_*\circ \phi_H\circ C_*^{H}(\alpha\times \beta) \\
&&= p_*\circ C_*^{G}\circ \phi_F(\alpha\times \beta)
= C_*^{G}\circ p_*\circ \phi_F(\alpha\times \beta) 
= C_*^{G}\circ\tau_F(\alpha, \beta).
\end{eqnarray*}
Thus the equality (\ref{pc}) is proved. This shows that $C^{sym}_*$ is 
a $\Q$-algebra homomorphism. 
Furthermore, $C^{sym}_*$ satisfies the naturality: 
$f^{sym}_*\circ C^{sym}_*=C^{sym}_*\circ f^{sym}_*$ 
for any proper morphism $f:X \to Y$, 
that immediately follows from 
$f^n_*\circ C^{S_n}_*=C^{S_n}_*\circ f^n_*$ for any $n$. 
Thus $C^{sym}_*$ is a natural transformation. 
\qed

\begin{rem} \label{integralcompatible}
{\rm 
(The degree) 
For $\alpha \in \F_{X,sym}[[z]]$, the integral 
$\int \alpha \in \Q[[z]]$ is equal to the $0$-th degree of $C^{sym}_0(\alpha)$, 
that is the value of pushforward $pt^{sym}_*C^{sym}_*(\alpha)$ 
induced by $pt: X \to \{pt\}$. 
}
\end{rem}

\begin{rem} \label{nonsingularcase} 
{\rm 
(Nonsingular case) 
Set 
$H^{X, sym}[[z]] :=\sum_{n=0}^\infty z^n H^{*}_{S_n}(X^n)$. 
The multiplication $\odot$ is given as 
$\omega \odot \omega':=
\frac{1}{|S_{m+n}|} \cdot \sum \sigma^*(p_1^*\omega \cup p_2^*\omega')$. 
where $p_1$ and $p_2$ are projections of $X^{m+n}=X^m\times X^n$ to the factors. 
For nonsingular $X$,  the Poincar\'e duality for all $n$ give 
${\cal P}: H^{X, sym}[[z]] \simeq H_{X, sym}[[z]]$ 
preserving the multiplications. 
Given a proper morphism $X \to Y$ between nonsingular varieties, 
the Gysin homomorphism 
$f_!^{sym}={\cal P}^{-1}\circ f^{sym}_*\circ {\cal P}$ is defined. 
}
\end{rem}

\subsection{Cycle types}  \label{cycletypes} 
A collection ${\bf c}=[c_1, \cdots, c_n]$ of non-negative integers 
satisfying $\sum_{i=1}^n ic_i=n$ 
is called {\it a cycle type of weight $n$}. 
We associate to a cycle type ${\bf c}$ 
three numbers, {\it weight}, {\it length} and {\it cardinality}, respectively, 
$$
|{\bf c}|:=\sum_{i=1}^n ic_i=n, \; \; 
l({\bf c}):=\sum_{i=1}^n c_i, \; \; 
\sharp {\bf c}:=\frac{n!}{1^{c_1}c_1!2^{c_2}c_2!\cdots n^{c_n}c_n!}.
$$
A representation $\rho: A \to S_n$ 
causes the orbit-decomposition of the $n$ points $\{1, 2, \cdots, n\}$. 
We say that $\rho$ is {\it of type ${\bf c}=[c_1, \cdots, c_n]$} 
if the decomposition consists of 
$c_k$ subsets (orbits) having exactly $k$ points. 
For ${\bf c}$ of weight $n$, 
let  ${\rm Hom}(A,S_n;{\bf c})$ denote 
the set of $S_n$-representations of type ${\bf c}$.  

For instance, ${\rm Hom}(\Z,S_n;{\bf c})$ is the set of 
permutations of type ${\bf c}$, i.e.  a conjugacy class in $S_n$; 
it consists  of $\sharp {\bf c}$ elements. 

Recall that 
$\jeden_{X^n/S_n}^{(A)}:=\frac{1}{n!}\sum \jeden_{(X^n)^{\rho(A)}}$ 
taken over all $S_n$-representations $\rho$ 
and $j_r(A)$ denotes the number of subgroups with index $r$.  
We denote by $\jeden_{{\mit\Delta} X^n}$ ($\in F^{S_n}_{inv}(X^n)$)  
the characteristic function of the diagonal ${\mit\Delta} X^n \subset X^n$. 
\begin{lem}\label{Dey}
The canonical  function 
$\jeden_{X^n/S_n}^{(A)}$ is equal to the sum
$$\sum_{|{\bf c}|=n}\; \frac{\sharp [c_1, \cdots, c_n]}{n!} \cdot  
 (j_1(A)\cdot \jeden_{{\mit\Delta}X})^{c_1}\odot \cdots \odot 
(j_n(A)\cdot \jeden_{{\mit\Delta}X^n})^{c_n},$$
taken over all cycle types of weight $n$. 
\end{lem}
\noindent{\it Proof}. 
This is elementary. 
Set $\theta_r:=\sum \jeden_{(X^r)^{\tau(A)}}$ 
taken over all transitive actions $\tau$ 
of $A$ on $r$ points $\{1, \cdots, r\}$. 
It is enough to show the following equalities 
(summing up $(2)$ over all $\bf c$ completes the proof): 
$$
\theta_r= (r-1)! j_r(A) \cdot \jeden_{\mit\Delta X^r}  \eqno{(1)}
$$
$$
\sum_{{\rm Hom}(A,S_n; {\bf c})} \jeden_{(X^n)^{\rho(A)}} 
= \frac{n!}{\prod_{r=1}^n (r!)^{c_r}c_r!} \cdot  
(\theta_1)^{c_1}\odot \cdots \odot (\theta_n)^{c_n}.
\eqno{(2)}
$$
It is well known (cf., \cite{Sta}) that 
the number of transitive actions of $A$ on $r$ points 
is given by $(r-1)! j_r(A)$ (in fact, 
given a subgroup $B$ with index $r$, 
an order of proper cosets of $A/B$, 
say $A=B \cup B_2 \cup \cdots \cup B_r$, determines  
a transitive action on $r$ points so that elements in $B_i$ send $1$ to $i$). 
It follows from the transitivity that 
$\jeden_{(X^r)^{\tau(A)}}= \jeden_{\mit\Delta X^r}$, thus $(1)$ is proved. 
To show $(2)$, note that 
every element $\rho$ in ${\rm Hom}(A, S_n;{\bf c})$ is uniquely obtained 
from a decomposition of $n$ points into 
disjoint subsets according to the type ${\bf c}$ 
and a choice of a transitive action on each of the subsets. 
Then, 
for such $\rho$,  
there are exactly $\prod (r!)^{c_r}c_r!(=:q)$ permutations $\sigma$ 
such that 
$\jeden_{(X^n)^{\rho(A)}}=
\sigma_*((\jeden_{\mit\Delta X^1})^{\times c_1}\times 
\cdots \times (\jeden_{\mit\Delta X^n})^{\times c_n})$. 
By definition we see 
$$
\sum_{\sigma \in S_n} 
\sigma_*((\theta_1)^{\times c_1}\times \cdots \times (\theta_n)^{\times c_n})
=n! \cdot 
(\theta_1)^{\odot c_1}\odot \cdots \odot (\theta_n)^{\odot c_n},
$$
which is just the sum of 
$\jjeden_{(X^n)^{\rho(A)}}$ over all $\rho$ of type ${\bf c}$ 
but taken account of $q$ times for each $\rho$, thus $(2)$ follows. 
\qed

\begin{rem}\label{conjugacy} 
{\rm 
Let $u_d(A)$ denote the number of conjugacy classes 
of subgroups in $A$ with index $d$. 
It holds {\rm (\cite{Sta}, Prob. 5.13)} that 
$$
j_k(A\times \Z)=\sum_{d|k} d \; u_d(A). 
$$
In particular, if $A$ is abelian, 
$u_d(A)=j_d(A)$ hence $j_k(A\times \Z)=\sum_{d|k} d \; j_d(A)$. 
Throughout let us write $j(m;k):=j_k(\Z^m)$ for short 
and set $j(0;1)=1$, $j(0;k)=0$ $(k >1)$. 
Then by induction, 
$j(m;k)=\sum j_1^{m-1}j_2^{m-2}\cdots j_{m-1}$
the sum taken over all  $m$-tuples of integers $(j_1, \cdots, j_m)$  
with $\prod_{i=1}^m j_i=k$. 
}
\end{rem}

\subsection{Diagonal operators}  \label{Diagonal} 
{\it The standard $n$-th diagonal operator}   
$D^n$ ($n=0, 1,\cdots$) is  defined as 
 the pushforward homomorphisms: 
$D^0=1$, $D^1=D=id_*$ ($id:X \to X$) and  $D^n:=({\mit\Delta}^n)_*$ 
$$
D^n: \F(X) \to \F^{S_n}_{inv}(X^n), 
\;\;  
D^n: H_*(X) \to H^{S_n}_*(X^n), 
$$ 
where ${\mit\Delta}^n: X \to X^n$  is the diagonal inclusion map, 
${\mit\Delta}^n(x):=(x, \cdots, x)$, regarded to be $S_n$-equivariant 
(with the trivial action on $X$).  
We call  a formal power series in $zD$, 
$U=\sum_{n=0}^\infty v_n z^nD^n$, with rational coefficients 
{\it a standard formal diagonal operator} . 
The algebra consisting of standard formal diagonal operators is denoted 
by $\Q[[zD]]$ (as a convention, $zD=Dz$).  
We put $\Q[[zD]]^+:=(zD)\Q[[zD]]$ (consisting of $U$ with zero constant term). 
There are canonically defined exponential and logarithmic functions: 
We denote them by 
$Exp(U)$ and $Log(1+U)$ for $U \in \Q[[zD]]^+$. 

 {\it The mixed $n$-th diagonal operator 
 of type ${\bf c}=[c_1, \cdots, c_n]$} means the maps 
$$
{\cal D}^{{\bf c}}: \F(X) \to \F^{S_n}_{inv}(X^n),  \quad 
{\cal D}^{{\bf c}}: H_*(X) \to H^{S_n}_*(X^n),
$$
given by 
$$
{\cal D}^{[c_1, \cdots, c_n]}(\alpha)
:=(D^1(\alpha))^{c_1}\odot (D^2(\alpha))^{c_2}\odot \cdots \odot
(D^n(\alpha))^{c_n},
$$ 
where 
$(D^k(\alpha))^c$ is the $c$ times multiple $(D^k)^{\odot c}$ on $\alpha$, 
that is, $D^k(\alpha) \odot \cdots \odot D^k(\alpha)$. 
We also define  {\it a formal diagonal operator}  as 
a formal series $T=\sum_{n=0}^\infty z^n T_n$ of  linear combinations 
$T_n=\sum_{|{\bf c}|=n} v_{\bf c} {\cal D}^{\bf c}$. 
All formal diagonal operators form a $\Q$-algebra with the multiplication $\odot$, 
which contains $\Q[[zD]]$ as a $\Q$-linear subspace. 
Every formal operator $T$ naturally induces maps 
$T:\F(X) \to \F_{X, sym}[[z]]$ and  $T: H_*(X) \to H_{X, sym}[[z]]$, 
which are linear if and only if $T$ is standard. 
By using $\odot$ we define 
$exp(T):=\sum_{n=0}^\infty \frac{1}{n!} T^{\odot n}$
for $T$ with zero constant term, as well $log(1+T)$. 
Obviously, $exp(T+T')=exp(T)\odot exp(T')$, etc. 
Since we use functions $exp$ and $Log$ only, 
the following notation would not cause any confusion: 

\

\noindent{\it Notation}: $\quad$  
$(1+U)^{\alpha}:=exp(Log(1+U)(\alpha))$ for $U \in \Q[[zD]]^+$. 
\begin{rem} \label{remarkT}
{\rm 
For instance, as formal operators, $(1-zD)^{-\alpha}\not= 
(1+zD+z^2D^2+\cdots)(\alpha)=Exp(-Log(1-zD))(\alpha)$. 
When $X=pt$, $\F_{pt,sym}[[z]]=H_{pt,sym}[[z]]=\Q[[z]]$, 
so $T$ defines $\Q \to \Q[[z]]$ $(z^n v {\cal D}^{\bf c}(a)=va^{l({\bf c})} z^{n})$. 
In this case, 
if $T$ is standard $(l({\bf c})=1)$ or $T=exp(U)$ for a standard $U$, 
the letter $D$ in $T$ has no longer meaning. For instance, $(1-zD)^a=(1-z)^a$. 
}
\end{rem}
\begin{prop} \label{Tpushforward} 
Let $T$ be a formal diagonal operator. Then, 
\begin{enumerate}
\item[(1)]
for a proper morphism $f:X \to Y$, it holds that 
$T\circ f_*=f^{sym}_*\circ T$ in both cases of constructible functions and homologies. 
In particular, $\int$ and $T$ commutes; 
\item[(2)]
 The following diagram commutes: 
$$
\begin{array}{ccc}
 \F(X) & \stackrel{C_*}{\longrightarrow} & H_*(X)  \\
T \;\; \downarrow \; \;\; \; & & \; \; \; \;\downarrow \; \; T \\
\F_{X, sym}[[z]] & \stackrel{C^{sym}_*}{\longrightarrow} & H_{X, sym}[[z]]
\end{array}
$$
\end{enumerate}
\end{prop}
\noindent{\it Proof}.   
Obviously, $D^n\circ f_*=(f^n)_*\circ D^n$, thus (1) holds. 
Since $S_n$ acts on $X$ trivially, 
$C^{S_n}_*$ for $X$ coincides with the ordinary $C_*:\F(X)\to H_*(X)$, and hence 
the naturality implies that $C^{S_n}_*\circ D^n=D^n\circ C^{S_n}_*=D^n\circ C_*$. 
Thus (2) follows.  
\qed

\begin{rem} \label{Diag}
{\rm 
We may denote as 
$T=\sum_{n=0}^\infty  T_n:\F(X)\otimes \Q[[z]] \to \F_{X,sym}[[z]]$ (as well homologies): 
for $\varphi(z)=\sum_{n=0}^\infty \alpha_nz^n$  
$(\alpha_n \in \F(X))$, 
$$
{\cal D}^{{\bf c}}(\varphi(z))
:=(D^1(\alpha_1))^{c_1}\odot \cdots \odot (D^n(\alpha_n))^{c_n}z^n.
$$
For example, $Log(1-zD)({\alpha})= 
- \sum_{r=1}^\infty \frac{z^r}{r}D^r(\alpha) 
=Log(1-D)(\alpha \otimes \sum_{r=1}^\infty z^r)$. 
As in Proposition \ref{Tpushforward}, 
the naturality and the compatibility are clear also in this sense. 
}
\end{rem}

\subsection{Generating functions}  
Now we are ready to prove Theorem \ref{DW} in Introduction. 
\begin{prop}\label{DWop} 
It holds that 
$$
\sum_{n=0}^\infty \jeden_{X^n/S_n}^{(A)} \, z^n
=exp\left( \sum_{r=1}^\infty \frac{j_r(A)}{r}z^rD^r (\jeden_X)\right).
$$
\end{prop}
\noindent{\it Proof}.   
A direct computation. 
For any standard $U=\sum_{k=1}^\infty \frac{a_k}{k}z^kD^k$, we see 
\begin{eqnarray*}
exp(U(\jeden_X))&=&
exp \left((a_1zD +\frac{a_2}{2}z^{2}D^{2}+ \cdots)\jeden_X \right) \\
&=& exp(a_1z \cdot \jeden_{\mit\Delta X}) \odot 
exp\left(\frac{a_2}{2}z^{2}\cdot \jeden_{\mit\Delta X^2}\right) \odot\cdots \\
&=&
\left(\sum_{c_1=0}^\infty 
\frac{a_1^{c_1}z^{c_1}}{1^{c_1}c_1!}(\jeden_{\mit\Delta X})^{c_1}\right) \odot 
\left(\sum_{c_2=0}^\infty 
\frac{a_2^{c_2}z^{2c_2}}{2^{c_2}c_2!}
(\jeden_{\mit\Delta X^2})^{c_2}\right) \odot \cdots \\
&=&
1+\sum_{n=1}^\infty z^{n}
\sum_{|{\bf c}|=n}
\frac{\sharp{\bf c}\;}{n!}\cdot  
(a_1\jeden_{\mit\Delta X})^{c_1}\odot\cdots \odot  
(a_n\jeden_{\mit\Delta X^n})^{c_n}. 
\end{eqnarray*}
Taking $a_i=j_i(A)$, then Lemma \ref{Dey} shows the assertion. 
\qed

\noindent{\it Proof of Theorem \ref{DW}}: 
Apply $C^{sym}_*$ to 
the both sides of this equality of constructible functions, then 
Proposition \ref{Tpushforward} (2) implies the theorem. 
\qed

\begin{exam} \label{freeabelian1} 
{\rm 
The generating function for $\jeden_{X^n/S_n}^{(A\times \Z)}$ 
is written as an exponential using the coefficients $u_d(A)$ 
instead of $j_r(A\times \Z)$ {\rm (Remark \ref{conjugacy})}.  
In case of $A=\Z^m$, 
 (letting $j(m;k)=j_k({\Z^m})$ and $Z=zD$) 
\begin{eqnarray*}
&&\sum_{k=1}^\infty \, j(m;k) \, \frac{Z^k}{k} (\alpha)
\, = \,  
\sum_{k=1}^\infty (\sum_{r|k} \,  j({m-1};r) r )\,\frac{Z^k}{k} (\alpha)  \\
&&\, = \,  
\sum_{r=1}^\infty \sum_{j=1}^\infty \,   j({m-1};r) \, \frac{Z^{rj}}{j} (\alpha) 
 \, = \,  
 \sum_{r=1}^\infty \, \, Log(1-Z^r) (- j({m-1};r) \alpha). 
\end{eqnarray*}
Thus, according to our convention of the notation, 
Proposition \ref{DWop} says 
$$
\sum_{n=0}^\infty \jeden_{X^n/S_n}^{(m)} \, z^n
=\prod_{r=1}^\infty  (1-z^rD^r)^{- j(m-1;r)\, \jjeden_X}, 
$$
and applying $C^{sym}_*$ to this gives 
$$
\sum_{n\ge 0}C_*^{S_n}(\jeden_{X^n/S_n}^{(m)}) \,  z^n 
= \prod_{r=1}^\infty 
(1-z^rD^r)^{-j(m-1;r)C_*(X)}. 
$$
The formulae (\ref{MacdonaldChern}) and  (\ref{HHChern}) 
in Introduction are the case of $m=1,2$, respectively. 
Let us see the $0$-dimensional homology part 
or equivalently the integral $\int$. 
Then by Proposition \ref{Tpushforward} $(1)$ (and also Remark \ref{remarkT} $(2)$) 
we recover Bryan-Fulman's formula \cite{BF}: 
$$
\sum_{n=0}^\infty 
\chi_m(X^n;S_n)\, z^n = 
\prod_{j_1, \cdots, j_{m-1}\ge 1} \, 
(1-z^{ j_1 j_2 
\cdots j_{m-1} })^{ - j_1^{m-2}j_2^{m-3}\cdots j_{m-2} \, \chi(X)}. 
$$ 
}
\end{exam}
\begin{exam} 
{\rm 
We note another typical cases, where $j_r(A)$ are easily counted. 

\begin{enumerate}
\item $A=\Z/d\Z$, the cyclic group of order $d$: 
$$
\sum_{n=0}^\infty C^{S_n}_*(\jeden_{X^n/S_n}^{(\Z/d\Z)}) \, z^n 
= exp\left(\sum_{r|d} \frac{1}{r}(zD)^r C_*(\jeden_X)\right). 
$$
\item $A=\Z_p$, the (additive) group of the $p$-adic integers ($p$ a prime): 
$$
\sum_{n=0}^\infty C^{S_n}_*(\jeden_{X^n/S_n}^{(\Z_p)}) \, z^n 
= exp\left(\sum_{k=0}^\infty \frac{1}{p^k}(zD)^{p^k} C_*(\jeden_X)\right).
$$
We may call the right-hand side 
``the Artin-Hesse exponential for the Chern class of $X$". 
\end{enumerate}
}
\end{exam}

\section{Wreath Products}

We discuss on a $G$-version of $C_*^{sym}$. 
Although the combinatorics looks much involved than the case of $S_n$, 
the construction is rather straightforward. 
Let $G$ be a finite group and $X$ a $G$-variety. 

\subsection{Action of wreath products} 
Let $G^n$ be the direct product of $n$ copies of $G$ 
and its element is denoted by $\bar{g}=(g_1, \cdots, g_n)$. 
The direct product of groups $G^n\times S_n$ acts on $X^n$ in an obvious way: 
$(\bar{g}, \sigma)$ sends $(x_1, \cdots)$ to $\sigma(g_1.x_1, \cdots)$ 
($=(g_{\sigma^{-1}(1)}.x_{\sigma^{-1}(1)}, \cdots)$). 
{\it The wreath product}, denoted by $G_n$ (or $G \sim S_n$), 
is defined as the semi-direct product 
of $G^n$ and $S_n$ with the multiplication 
$(\bar{h}, \sigma)(\bar{g}, \tau):=(\bar{h} \cdot\sigma(\bar{g}), \sigma\tau)$, 
where 
$\sigma(\bar{g}):=\bar{g}\circ \sigma^{-1}=(g_{\sigma^{-1}(1)}, \cdots)$. 
The action of $G_n$ on $X^n$ is given as 
$$
(g_1, \cdots, g_n,\sigma).(x_1,\cdots, x_n):=
(g_1.x_{\sigma^{-1}(1)}, \cdots, g_n.x_{\sigma^{-1}(n)}).
$$ 

A $G_n$-representation $\rho: A \to G_n$ is denoted by $(\bar{\rho}, \sigma)$ 
where $\bar{\rho}(a) \in G^n$, $\sigma(a) \in S_n$ for $a\in A$. 
We say that $\rho$ is {\it of type ${\bf c}=[c_1, \cdots ,c_n]$} 
if the second factor $\sigma$ 
gives a $S_n$-representation of cycle type ${\bf c}$. 
For any ${\bf c}$ of weight $n$, 
we denote by ${\rm Hom}(A, G_n;{\bf c})$ the set of 
$G_n$-representations of $A$ of type ${\bf c}$.

\subsection{Formal power series}
At first we set 
$$
\F^G_{X, sym}[[z]]:=\sum_{n=0}^\infty z^n \F^{G_n}_{inv}(X^n), \quad 
H^G_{X, sym}[[z]] :=\sum_{n=0}^\infty z^n H^{G_n}_{2*}(X^n).
$$
It is easy to see that they become commutative and associative graded $\Q$-algebras 
by mean of the multiplication 
$$
\alpha \odot_G \beta
:=\frac{1}{|G|^{m+n} (m+n)!} \sum_{(\bar{g}, \sigma) \in G_{m+n}} 
(\bar{g}, \sigma) _*(\alpha \times \beta)
$$
for $\alpha \in \F^{G_m}_{inv}(X^m)$ and 
$\beta \in \F^{G_{n}}_{inv}(X^{n})$ (as well the homology case). 
For a proper $G$-morphism $f:X\to Y$, 
the pushforward $f^{G, sym}_*$  is immediately defined. 
\begin{thm} 
$C^{G,sym}_*
=\prod_{n=0}^\infty C^{G_n}_*:\F^G_{X, sym}[[z]] \to H^G_{X, sym}[[z]]$ 
is a natural transformation. 
If the $G$-action on $X$ is trivial, then $C^{G,sym}_*$ coincides with $C^{sym}_*$. 
\end{thm}
\noindent{\it Proof}. 
 It is shown in the entirely same way as in the proof of Theorem \ref{SymC}. 
 \qed

\subsection{Diagonal operators}
{\it The standard operator}  
$D^n: \F_{inv}^G(X) \to \F_{inv}^{G_n}(X^n)$  is defined by 
$$
D^n(\alpha):=
\frac{1}{|G|^nn!} 
\sum_{(\bar{g},\sigma) \in G_n} (\bar{g}, \sigma)_*(\mit\Delta^n)_*(\alpha), 
$$
as well the homology case. If $G=\{e\}$, this $D^n$ coincides 
with the previous one in \ref{Diagonal}. 
Obviously, $D^n(\alpha)=1/|G|^n \cdot \sum_{\bar{g}} \bar{g}_*(\mit\Delta^n)_*(\alpha)$, 
being also invariant under the action of the direct product $G^n\times S_n$. 

Using the above $D^n$ and $\odot_G$, 
{\it formal diagonal operators} are defined in the same way as before: 
$T=\sum T_n$, $T_n=\sum v_{\bf c}{\cal D}^{\bf c}$ 
(see Remark \ref{Diag}). 
These formal operators $T$ enjoy 
the same properties as in Proposition \ref{Tpushforward} 
(the proof is the same). 
In particular, the following diagram commutes: 
$$
\begin{array}{rcl}
\F_{inv}^G(X)\otimes \Q[[z]] & \stackrel{C^G_*\otimes id}{\longrightarrow} 
& H_*^{G}(X) \otimes \Q[[z]]\\
T\downarrow \qquad \quad & & \;\; \qquad  \downarrow T \\
\F^G_{X,sym}[[z]] \quad & 
\stackrel{C^{G,sym}_*}{\longrightarrow} & \quad  H^G_{X,sym}[[z]]
\end{array}
\eqno{(d)}
$$

Recall the definition of canonical constructible functions, 
$$
\jeden_{X^n/G_n}^{(A)} := 
\frac{1}{|G|^n n!} \sum_{{\rm Hom}(A,G_n)} \jeden_{(X^n)^{\rho(A)}}.
$$ 
We have the following formula and the proof is given in the next subsection: 
\begin{prop}\label{DWGop} 
Under the same assumption as in Theorem \ref{DWG}, 
it holds that 
$$
\sum_{n=0}^\infty \jeden_{X^n/G_n}^{(A)} \, z^n
=exp\left( \sum_{r=1}^\infty \sum_{B \in \Omega_A(r)} 
\frac{1}{r}z^rD^r (\jeden_{X/G}^{(B)})\right).
$$
\end{prop}
In case of $G=\{e\}$ $\jeden_{X/G}^{(B)}=\jeden_X$,  
so this is the same as Proposition \ref{DWop}. 

\

 \noindent{\it Proof of Theorem \ref{DWG}}
It immediately follows from Proposition \ref{DWGop} 
and the commutative diagram $(d)$.  
\qed

\begin{rem} \label{Muller}
{\rm 
When $X=pt$, the above proposition recovers the enumerative formula 
$$
\sum_{n=0}^\infty \frac{|{\rm Hom}(A,G_n)|}{|G|^n n!} \, z^n 
= exp\left(\sum_{B \in \Omega_A} 
\frac{|{\rm Hom}(B,G)|}{|G| \cdot |A:B|}z^{|A:B|}\right)
$$
which has been given by M\"uller \cite{Mu} (Cor. 1, Ex. 4). 
If $G=\{e\}$,  it is the classical Dey-Wohlfahrt formula \cite{W}. 
}
\end{rem}

\begin{exam}\label{freeabelian2}
{\rm 
If $A=\Z^m$, then any subgroup $B$ of finite index turns again 
a lattice of rank $m$; therefore $\jeden_{X/G}^{(B)}=\jeden_{X/G}^{(m)}$. 
Thus it follows from Proposition \ref{DWGop} that 
$$
\sum_{n=0}^\infty \jeden_{X^n/G_n}^{(m)} \, z^n
=\prod_{r=1}^\infty  (1-z^rD^r)^{- j(m-1;r)\, \jjeden_{X/G}^{(m)}}. 
$$
Then $C^{G, sym}_*$  provides the Chern class formula; in particular, 
we recover Tamanoi's formula  \cite{Ta}: 
$$
\sum_{n=0}^\infty \chi_m(X^n; G_n) \, z^n
=\prod_{r=1}^\infty  (1-z^r)^{- j(m-1;r)\,  \chi_m(X; G)}. 
$$
As a simple comparison, 
the direct product case  essentially goes to Proposition \ref{DWop} 
and the result differs from the above for $m \ge 2$: 
$$
\sum_{n=0}^\infty \jeden_{X^n/G^n\times S_n}^{(m)} \, z^n
=\prod_{r=1}^\infty  (1-z^rD^r)^{- j(m-1;r)\, \jjeden_{X/G}^{(1)}}. 
$$
}
\end{exam}

\subsection{Canonical functions} \label{Gconst}
For our convenience, we put 
$$
\jeden_{X/G}^{(A;r)}:=\sum_{B\in \Omega_A(r)} \jeden_{X/G}^{(B)} 
\;  =  \; \frac{1}{|G|}\sum_B\sum_{\nu} \jeden_{X^{\nu(B)}} 
$$
where $\nu$ runs all $G$-representations of $B$. 
Then we have  
\begin{lem} \label{DeyG}
It holds that 
$$\jeden_{X^n/G_n}^{(A)} 
= \sum_{|{\bf c}|=n} \frac{\sharp{\bf c}\;}{n!} 
\left(D^1(\jeden_{X/G}^{(A;1)})\right)^{c_1} \odot_G \cdots \odot_G 
\left(D^n(\jeden_{X/G}^{(A;n)})\right)^{c_n}.$$ 
\end{lem} 
This lemma generalizes Lemma \ref{Dey}: 
By definition, if $G=\{e\}$, then 
$$
D^r(\jeden_{X/G}^{(A;r)})=D^r(j_r(A)\cdot \jeden_X)
=j_r(A)\cdot \jeden_{\mit\Delta X^r}.
$$
\noindent{\it Proof}.   
Set 
$\Theta_r:=\sum \jeden_{(X^r)^{\tau(A)}}$ taken over all 
$\tau \in {\rm Hom}(A,G_r;[0,\cdots, 0,1])$, 
that is, $G_r$-representations with transitive permutations.  
It suffices to show that 
$$
\frac{1}{|G|^r} \Theta_r
=(r-1)! \cdot D^r(\jeden_{X/G}^{(A;r)}) \eqno{(1')}
$$
$$
\sum_{{\rm Hom}(A,G_n;{\bf c})} \jeden_{(X^n)^{\rho(A)}}
= \frac{n!}{\prod_{r=1}^n (r!)^{c_r}c_r!}\cdot 
(\Theta_1)^{c_1} \odot_G \cdots \odot_G (\Theta_n)^{c_n}.
\eqno{(2')}
$$
The same proof as seen before (Lemma \ref{Dey} $(2)$) 
works also for $(2')$. So we now show $(1')$. 
Since the both sides of $(1')$ are $G_r$-invariant, 
it is enough to show that their values coincide 
at $\mbox{\boldmath $x$}=(x, \cdots,x) \in \mit\Delta X^r$. 
We denote the stabilizer group of $x$ of the $G$-action by $H:=Stab_G(x)$ 
and the $G$-orbit of $x$ by $G.x$. 
The RHS of $(1')$ at $\mbox{\boldmath $x$}$ is written as 
$$
(r-1)! \cdot D^r(\jeden_{X/G}^{(A;r)}) (\mbox{\boldmath $x$}) 
= \frac{(r-1)!\cdot |H|^{r-1}}{|G|^r}\sum_B |{\rm Hom}(B,H)|. 
$$
In fact, the number of $\bar{g} \in G^r$ so that $g_1.x=\cdots = g_r.x$ 
is equal to $|H|^r\cdot |G.x|$, so 
\begin{eqnarray*}
D^r(\jeden_{X/G}^{(A;r)}) (\mbox{\boldmath $x$}) 
&=&\frac{1}{|G|^{r+1}}
\sum_{B,\, \nu} \sum_{\bar{g}}
\bar{g}_*(\jeden_{\mit\Delta (X^{\nu(B)})^r})(\mbox{\boldmath $x$}) \\
&=& 
\frac{1}{|G|^{r+1}}\cdot |H|^r\cdot |G.x| \cdot
\sum_{B, \, \nu} \jeden_{X^{\nu(B)}}(x). 
\end{eqnarray*}

On one hand, the LHS of $(1')$ at $\mbox{\boldmath $x$}$ means 
$$
\frac{1}{|G|^r} \Theta_r(\mbox{\boldmath $x$}) 
= \frac{1}{|G|^r} |{\rm Hom}(A,H_r; [0, \cdots, 0,1])|, 
$$
where $H_r$ is the wreath product $H\sim S_r$. 
Consider the following composed map 
$$
{\rm Hom}(A,H_r; [0, \cdots, 0,1]) \to {\rm Hom}(A,S_r; [0, \cdots, 0,1]) 
\to \Omega_A(r)
$$
sending $(\bar{\rho}, \sigma) \mapsto \sigma \mapsto Stab_\sigma(1)$ 
(the stabilizer subgroup $\{a \in A, \sigma(a)(1)=1\}$). 
As seen in the proof of Lemma \ref{Dey}, 
for  each subgroup $B$ with $|A:B|=r$,  there are $(r-1)!$ choices of $\sigma$ mapped to it. 
Now fix such a $S_r$-representation $\sigma$. 
We take $a_i \in A$  ($2\le i \le r$) so that $\sigma(a_i)(1)=i$, 
and put 
$a_1=e'$ and $g_1=e$, the identities of $A$ and $H$ respectively. 
 Then any $\rho=(\rho_1, \cdots, \rho_r, \sigma): A \to H_r$ (a lift of $\sigma$) 
is uniquely determined by a choise of a representation $\nu: B \to H$ 
and ordered $(r-1)$ elements $g_i \in H$ ($2\le i \le r$).  
In fact,  a simple computation shows that 
$\rho_k$ ($1 \le k \le r$) has the following unique form  
with respect to  $\sigma$ and $\{a_i\}$: for $a \in A$, 
$\rho_k(a)=g_k\nu(a_k^{-1}aa_i)g_{i}^{-1} \in H$ 
where $i=\sigma(a)^{-1}(k)$. 
Hence 
$$
|{\rm Hom}(A,H_r; [0, \cdots, 0,1])| 
= \sum_B (r-1)! \cdot |H|^{r-1} \cdot |{\rm Hom}(B,H)|, 
$$
thus $(1')$ holds. This completes the proof. 
\qed

\noindent{\it Proof of Proposition \ref{DWGop}}. 
It suffices to expand 
$$
exp(-Log(1-D)\varphi_A(z)) \;\; \mbox{where} \;\; 
\varphi_A(z)=\sum_{r=1}^\infty \jeden_{X/G}^{(A;r)} z^r 
$$
like as the proof of Proposition \ref{DWop} (see Remark \ref{Diag}).  \qed

\end{document}